\documentclass[11pt,english]{article}
\usepackage[T1]{fontenc}
\usepackage[latin9]{inputenc}
\usepackage{geometry}
\usepackage{verbatim}
\usepackage{float}
\usepackage{rotfloat}
\usepackage{amsmath}
\usepackage{amssymb}
\usepackage{graphicx}
\usepackage{esint}

\makeatletter

\usepackage{times}

\newtheorem{proposition}{}

\newcommand{\ackname}{Acknowledgements:}

\usepackage{babel}

\makeatother

\usepackage{babel}
\begin{document}

\title{ON ILL-POSEDNESS OF NONPARAMETRIC INSTRUMENTAL VARIABLE REGRESSION  \\
 WITH CONVEXITY CONSTRAINTS}

\author{Olivier Scaillet$^{a}${*}}

\date{This version: July 2016.}
\maketitle
\begin{abstract}
{\footnotesize{}{}{} This note shows that adding monotonicity or convexity constraints on the regression function does not restore well-posedness in nonparametric instrumental variable regression.
The minimum distance problem without regularisation is still locally ill-posed. \medskip{}
 }{\footnotesize \par}

\noindent \textit{\footnotesize{}{}{}JEL Classification:}{\footnotesize{}{}{}{}
C13, C14, C26.}{\footnotesize \par}

\noindent \textit{\footnotesize{}{}{}Keywords:}{\footnotesize{}{}{}{}
Nonparametric Estimation, Instrumental Variable, Ill-Posed Inverse Problems.}{\footnotesize \par}

\medskip{}

\noindent {\scriptsize{}{}{}$^{a}$University
of Geneva and Swiss Finance Institute.}{\scriptsize \par}

\noindent {\tiny $^{*}$ Acknowledgements:
We thank Xiaohong Chen and Daniel Wilhelm  for helpful comments.}{\scriptsize{}{}{}{} }{\scriptsize \par}
\end{abstract}
\pagebreak{}

\section{Introduction}

We consider estimation of the regression model $Y = \varphi_0 (X) + U$.  
The variable $X$  has compact support ${\cal X} = [0,1]$ and is potentially endogenous.  
The instrument $Z$ has compact support ${\cal Z} = [0,1]$. The parameter of interest is the function $\varphi_0$ defined on ${\cal X}$
which satisfies the nonparametric instrumental   variable regression (NPIVR):
\begin{eqnarray} 
E [Y - \varphi_0 (X) \vert Z ] = 0. \label{NIVR}
\end{eqnarray}
As shown in Example 1 of Gagliardini and Scaillet (2016),  we do not need independence between the error $U$ and the instrument $Z$. 
This exemplifies a difference between restrictions induced by a parametric conditional moment setting and their nonparametric counterpart.
NPIVR estimation has received considerable attention in the recent years
building on a series of fundamental papers on ill-posed endogenous mean regressions 
(Ai and Chen (2003),  Newey and Powell (2003), Hall and Horowitz (2005), 
Blundell, Chen, and Kristensen (2007), Darolles, Fan, Florens, and
Renault (2011), Horowitz (2011)), and the
review paper by Carrasco, Florens, and Renault (2007). The
main issue in nonparametric estimation with endogeneity is
overcoming ill-posedness of the associated inverse problem. It
occurs since the mapping of the reduced form parameter (that is,
the distribution of the data) into the structural parameter (that
is, the instrumental regression function) is not continuous in the conditional moment $E [Y \vert Z ]$. We
need a regularization of the estimation to recover consistency.
Gagliardini and Scaillet (GS, 2012a) study a
Tikhonov Regularized (TiR) estimator (Tikhonov (1963a,b), Groetsch
(1984), Kress (1999)). They achieve regularization by adding a
compactness-inducing penalty term, the Sobolev norm, to a
functional minimum distance criterion. Chen and Pouzo (2012) study nonparametric estimation of conditional moment restrictions in which the generalized residual functions can be nonsmooth in the unknown functions of endogenous variables. For such a nonparametric nonlinear instrumental variables problem, they propose a class of penalized sieve minimum distance estimators (see Chen and Pouzo (2015) for inference in such a setting).
As discussed in Matzkin (1994), in nonparametric models, we can use economic restrictions, as in parametric models, to reduce the variance of estimators, to falsify theories, and to extrapolate beyond the support of the data. But, in addition, we can use some economic restrictions to guarantee the identification of some nonparametric models and the consistency of some nonparametric estimators.
Economic theory often provides shape restrictions on functions of interest in applications, such as
monotonicity, convexity, non-increasing (non-decreasing) returns to scale but economic theory does not provide finite-dimensional parametric models. This
motivates nonparametric estimation under shape restrictions. Since nonparametric estimates are often
noisy,  shape restrictions helps to stabilize nonparametric estimates without imposing arbitrary restrictions (see the recent works of Blundell, Horowitz, and Parey (2012), Horowitz and Lee (2015)).
Following that line of thought, we could hope that adding monotonicity or convexity constraints on the regression function would help
 to restore well-posedness in nonparametric instrumental variable regression. 
 The next section shows that this is unfortunately not the case since the minimum distance problem without regularisation is still locally ill-posed. Chetverikov and Wilhelm (2015) look at imposing two monotonicity conditions: (i) monotonicity of the regression function $\varphi_0$ and (ii) monotonicity of the reduced form relationship between the endogenous regressor $X$ and the instrument $Z$ in the sense that the conditional distribution of $X$ given $Z$ corresponding to higher values of $Z$ first-order stochastically dominates the same conditional distribution corresponding to lower values of $Z$ (the monotone IV assumption).
They show that these two monotonicity conditions together significantly change the structure of the NPIV model, and weaken its ill-posedness. In particular they point out that, even if well-posedness is not restored, those two monotonicity constraints improve the rate of convergence in shrinking neighborhoods of the constraint boundary and can have a significant impact on the estimator finite sample behavior. Chen and Christensen (2013) show that imposing shape restrictions only is not enough to improve convergence rates as long as the derivative constraints hold with strict inequality (i.e., in the interior of the constraint space). There may be rate improvements when the constraint is binding.

\section{Ill-posedness with convexity constraints}

The functional parameter $\varphi_0$ belongs to a subset $\Theta$ of  $L^2({\cal X})$, where $L^{2}(\mathcal{X})$ denotes the $L^{2}$ space of square
integrable functions of $X$ defined by the scalar product $\langle
\varphi ,\psi \rangle =\int \varphi (x)\psi (x)dx$, and we write $\Vert \varphi \Vert$ for the $L^2$ norm $\langle
\varphi ,\varphi \rangle^{1/2}$.

We assume the following identification condition. \bigskip \newline
\textbf{Assumption 1: }$\varphi _{0}$\textit{\ is the unique function }%
$\varphi \in L^2({\cal X}) $\textit{\ that
satisfies the conditional moment restriction (\ref{NIVR}).}\bigskip \newline
We refer to Newey and Powell (2003), Theorems 2.2-2.4, for sufficient conditions ensuring
Assumption 1. 

Let us  consider a nonparametric
minimum distance approach to obtain $\varphi _{0}$. This relies on $\varphi_{0}$ minimizing 
\begin{equation}
Q_{\infty}(\varphi):=E\left[ m \left( \varphi, Z \right)^{2}\right],\, \, \varphi \in L^2({\cal X}),  \label{Qinfty}
\end{equation}%
where  $m \left( \varphi, Z \right) = E [Y - \varphi (X) \vert Z ]$.
We can write the
conditional moment function $m\left( \varphi ,z\right) $
 as: 
\begin{equation}
m\left( \varphi ,z\right) =\left( A \varphi \right) \left( z \right) -r \left( z\right) =\left(
A \Delta \varphi \right) \left( z \right) ,  \label{mphiZ}
\end{equation}%
with $\Delta \varphi :=\varphi -\varphi _{0},$ and  where the linear
operator $A$ is defined by $\left( A \varphi \right)
(z):=\break \displaystyle\int \varphi
(x)f_{X|Z}(x \vert z)dx$ and $r\left( z\right)
:=\displaystyle\int y f_{Y|Z}(y|z)dy,$ where $f_{X|Z}$ and $f_{Y|Z}$ are
the conditional densities of $X$ given $Z$, and $Y$ given $Z$.
Assumption 1 on identification of $\varphi _{0}$ holds if and only if
operator $A$ is injective.
Further, we assume that $A$ is a bounded operator from $L^{2}({\cal X}) $ to $L^{2}(\mathcal{Z})$, where $L^{2}(\mathcal{Z})$ denotes the $L^{2}$ space of square
integrable functions of $Z$ defined by the scalar product $\langle \psi
_{1},\psi _{2}\rangle _{L^{2}(\mathcal{Z})}=E\left[ \psi _{1}\left( Z\right) \psi _{2}\left( Z\right)\right] .$ The limit criterion (\ref{Qinfty}) becomes 
\begin{eqnarray}
Q_{\infty }(\varphi ) &=&\langle A\Delta \varphi,A \Delta \varphi \rangle _{L^{2}(\mathcal{Z})}  \label{limit2},\end{eqnarray}

\textbf{Assumption 2:} \textit{The linear operator }$A$ \textit{from 
}$L^{2}({\cal X}) $\textit{\ to }$L^{2}(\mathcal{Z%
})$\textit{\ is compact.\bigskip }\newline
Assumption 2 on compactness \ of operator $A$ holds under mild
conditions  on the conditional density $f_{X|Z}$ (see e.g. GS). In the proof of Proposition 1 below, we also need the regularity conditions: $\sup_z \vert f_Z (z) \vert < \infty$ and $\sup_{x,z} \vert f_{X\vert Z} (x\vert z) \vert < \infty$.

Proposition 1 shows that the minimum distance problem above is
locally ill-posed (see e.g.\ Definition 1.1 in Hofmann and Scherzer (1998)) even if we consider
monotonicity, monotonicity nonnegativity, or convexity constraints.
There are sequences of increasingly oscillatory functions arbitrarily close
to $\varphi _{0}$ that approximately minimize $Q_{\infty }$ while not
converging to $\varphi _{0}$. In other words, function $\varphi _{0}$ is not
identified in $\Theta $ as an isolated minimum of $Q_{\infty }$. Therefore,
ill-posedness can lead to inconsistency of the naive analog estimators based
on the empirical analog of $Q_{\infty }$. In order to rule out these
explosive solutions, we can use penalization as in GS (see Gagliardini and Scaillet (2012b) and Chen and Pouzo (2012) for the quantile regression case).
Under a stronger assumption than Assumption 1, namely local injectivity of $A$, the definition of local ill-posedness is
equivalent to $A^{-1}$ being discontinuous in a neighborhood of $A(\varphi_0)$ (see Engl,
Hanke, and Neubauer (2000, Chapter 10)).

\begin{proposition} Let $\varphi_0$ satisfy  monotonicity, monotonicity nonnegativity, or convexity constraints. Then, under Assumptions 1 and 2, the minimum distance problem is locally ill-posed,
namely for any $r > 0$ small enough, there exist $\varepsilon \in (0,r)$ and a sequence 
$(\varphi_n) \subset B_r(\varphi_0) := \{\varphi \in L^2({\cal X}) : \Vert \varphi - \varphi_0 \Vert < r \}$ such that  $\Vert \varphi_n -\varphi_0 \Vert \geq \varepsilon$,  $Q_{\infty} (\varphi_n) \rightarrow Q_{\infty}(\varphi_0) = 0$, and such that $\varphi_n$ satisfies the same constraints as $\varphi_0$.
\end{proposition}

\textbf{Proof:} The proof of Proposition 1 gives explicit sequences $(\varphi_n)$ generating ill-posedness when $\varphi_0$ satisfies  monotonicity, monotonicity nonnegativity, or convexity constraints.

Let us build $\varphi _{n}=\varphi
_{0}+\varepsilon \psi _{n}$, $\varepsilon
>0$, where $\psi _{n}(x) := - (2n+1)^{1/2} (1-x)^n$ and $\varphi_0$ is monotone and increasing.
Then $\varphi _{n} \in L^{2}({\cal X})$ and the first detivative $\nabla \varphi _{n} \geq 0$.
Since $\left\Vert \psi _{n}\right\Vert =1$, when we choose $\varepsilon
>0$ sufficiently small, we have $\left( \varphi _{n}\right) \subset
B_{r}(\varphi _{0})$, and\textit{\ }$\varphi _{n}\nrightarrow \varphi _{0}$. We also have that $A \varphi_n \stackrel{w}{\rightarrow} A \varphi_0$, where  $\stackrel{w}{\rightarrow}$ denotes weak convergence. Indeed,  for $q \in L^{2}(\mathcal{Z})$, we get $ \langle q ,A \varphi_n \rangle _{L^{2}(\mathcal{Z})} = \langle q ,A \varphi_0 \rangle _{L^{2}(\mathcal{Z})} + \varepsilon \langle q ,A \psi_n \rangle _{L^{2}(\mathcal{Z})}$, and $ \langle q ,A \psi_n \rangle _{L^{2}(\mathcal{Z})} \rightarrow 0$, since $\displaystyle \vert A \psi_n\vert \leq C (2n+1)^{1/2} \frac{1}{n+1}$ for $C >0$.
Since $A$ is
compact and $\left( \varphi _{n}\right) $ is bounded, the
sequence $A \varphi _{n} $ admits a convergent
subsequence $A \varphi _{m(n)} \rightarrow \xi $.
Since the weak limit is unique, we have $\xi =A \varphi
_{0}$. Thus $A\varphi _{m(n)}\rightarrow 
A \varphi _{0} $ and $Q_{\infty }\left( \varphi
_{m(n)}\right) \rightarrow 0$ but $\left\Vert \varphi _{m(n)}-\varphi
_{0}\right\Vert \geq \varepsilon $, hence the stated result follows.

The above argument works also with the function $\psi _{n}(x) := \displaystyle \frac{(2n+1)^{1/2}}{(2^{2n+1}-1)^{1/2}} (1+x)^n$,
which yields a monotone nonnegative function $\varphi _{n} \in L^{2}({\cal X})$ if $\varphi _{0} \geq 0$
and is monotone. 
This shows that the positivity constraint does not help 
here either. 

Since the higher order derivatives $\nabla^m \psi_{n} \geq 0$, $m \geq 1$, this example also shows that positivity constraints 
on the higher order derivatives $\nabla^m \varphi_0 \geq 0$, such as a convexity constraint  $\nabla^2 \varphi_0 \geq 0$, 
does not restore well-posedness of the estimation problem 
in the NPIVR setting.

\end{document}